\documentstyle[12pt]{article}
\newtheorem{Theorem}{Theorem}
\newtheorem{Lemma}{Lemma}

\makeatletter
\let\normalequation=\equation
\def\equation{\@ifnextchar[{\subequation}{\normalequation}}
\def\subequation[#1]#2{\@ifundefined{r@#1}%
  {\def\theequation{\bf ??#2}\@warning
    {Reference `#1' on page \thepage \space
     undefined}}{\edef\@tempa{\@nameuse{r@#1}}%
    \edef\theequation{\expandafter\@car\@tempa \@nil#2}}%
  \let\@currentlabel\theequation $$}
\makeatother
\begin{document}

\begin{center}
{\large\bf Global Exponential Stability of Almost Periodic
Solution for A Large Class of Delayed Dynamical
Systems}\footnote{It is supported
  by National Science Foundation of China 60374018 and 60074005.}
\end{center}

\begin{center}
Wenlian Lu \footnote{This author is with Lab. of Nonlinear
Mathematics Science, Institute of Mathematics, Fudan University,
Shanghai, 200433, P.R.China.} and Tianping Chen,
\footnote{Corresponding author, he is  with Lab. of Nonlinear
Mathematics Science, Institute of Mathematics, Fudan University,
Shanghai, 200433, P.R.China. Email: tchen@fudan.edu.cn}
\end{center}

\thispagestyle{empty}

\vspace{0.5in}

\begin{abstract}

Research of delayed neural networks with variable
self-inhibitions, inter-connection weights, and inputs is an
important issue. 
In this paper, we discuss a large class of delayed dynamical
systems with almost periodic self-inhibitions, inter-connection
weights, and inputs. This model is universal and includes delayed
systems with time-varying delays, distributed delays as well as
combination of both. We prove that under some mild conditions, the
system has a unique almost periodic solution, which is globally
exponentially stable. We propose a new approach, which is
independent of existing theory concerning with existence of almost
periodic solution for dynamical systems.

\end{abstract}

\section{\bf Introduction}

Recurrently connected neural networks, sometimes called
Grossberg-Hopfield neural networks, have been extensively studied
in past years and found many applications in different areas.
However, many applications heavily depend on the dynamic behaviors
of the networks. Therefore, analysis of these dynamic behaviors is
a necessary step toward practical design of these neural networks.

Recurrently connected neural network is described by the following
differential equations:
\begin{equation}
\frac{du_{i}(t)}{dt}=-d_{i}u_{i}(t)+\sum_{j}a_{ij}g_{j}(u_{j}(t))+I_{i}
\quad i=1,\cdots,n \label{Hopfield}
\end{equation}
where $g_{j}(x)$ are activation functions, $d_{i}$, $a_{ij}$ are
constants and $I_{i}$ are constant inputs. In practice, however,
the interconnections are  generally asynchronous. Therefore, one
often needs to investigate the following delayed dynamical
systems:
\begin{eqnarray}
\frac{du_{i}(t)}{dt}=-d_{i}u_{i}(t)+\sum_{j=1}^{n}a_{ij}g_{j}(u_{j}(t))
+\sum_{j=1}^{n}b_{ij}f_{j}(u_{j}(t-\tau_{ij}(t)))+I_{i} \quad
i=1,\cdots,n \label{delay}
\end{eqnarray}
where activation functions $g_{j}$ and $f_{j}$ satisfy certain
defining conditions,  and $a_{ij}$, $b_{ij}$, $I_{i}$ are
constants.

However, the interconnection weights $a_{ij}$, $b_{ij}$,
self-inhibition $d_{i}$ and inputs $I_{i}$ should be variable with
time.  Therefore, we need to discuss following dynamical systems.
\begin{eqnarray}
\frac{du_i}{dt}=-d_{i}(t)u_{i}(t)+\sum_{j=1}^{n}a_{ij}(t)g_j(u_j(t))
+\sum_{j=1}^{n}b_{ij}(t)f_{j}(u_{j}(t-\tau_{ij}(t)))+I_i(t)
\label{al}
\end{eqnarray}
or delayed systems with distributed delays and time-varying delays
\begin{eqnarray}
&&\frac{du_{i}(t)}{dt}=-d_{i}(t)u_{i}(t)+\sum_{j=1}^{n}a_{ij}(t)g_{j}(u_{j}(t))\nonumber\\
&+&\sum_{j=1}^{n}b_{ij}(t)\int_{0}^{\infty}k_{ij}(s)f_{j}(u_{j}(t-\tau_{ij}(t)-s))ds
+I_{i}(t) \quad i=1,\cdots,n
\end{eqnarray}
To unify these delayed systems, in \cite{Chen}, we proposed the
following  general dynamical systems
\begin{eqnarray}
&&\frac{du_i}{dt}=-d_{i}(t)u_{i}(t)+\sum_{j=1}^{n}a_{ij}(t)g_j(u_j(t))\nonumber\\
&+&\sum_{j=1}^{n}\int_{0}^{\infty}f_{j}(u_{j}(t-\tau_{ij}(t)-s))dK_{ij}(t,s)+I_i(t)\quad
i=1,\cdots,n\label{almostperiodic}
 \label{periodic}
\end{eqnarray}
where $\tau_{ij}\ge 0$ and $s\rightarrow dK_{ij}(t,s)$ is a
Lebesgue-Stieljies measure for each $t\in R$.

The initial conditions are
\begin{equation}
u_{i}(s)=\phi_{i}(s) \quad  for \quad s\in(-\infty,0]\label{ini}
\end{equation}
where $\phi_{i}\in C((-\infty,0]), i=1,\cdots,n$.

There are several papers discussing periodic dynamical systems and
their periodic solutions and its stability, see
\cite{Cao1,Chen,Mo,Go,Xiang,Zhang,Zheng,Zhou} and others. For
example, in \cite{Zheng,Zhou}, under assumptions that
$d_{i}(t)>d_{i}>0$, $a_{ij}(t), b_{ij}(t), I_i(t):
\mathbf{R}^{+}\rightarrow \mathbf{R}$ are continuously periodic
functions with period $\omega>0$, i.e.,
$d_{i}(t+\omega)=d_{i}(t)$, $a_{ij}(t)=a_{ij}(t+\omega)$,
$b_{ij}(t)=b_{ij}(t+\omega)$, $I_{i}(t)=I_{i}(t+\omega)$,
$\tau_{ij}(t+\omega)=\tau_{ij}(t)$ for all $t>0$ and
$i,j=1,2,\ldots,n$, the model (\ref{delay}) was investigated.
Existence of periodic solution and its stability were
investigated. And in \cite{Zheng}, comparisons of various
stability criteria were given  too. In \cite{Chen,Lu}, a general
approach to discuss existence of periodic solution and its
stability analysis for the model (\ref{al}) was proposed. In
\cite{Huang,A1}, authors presented some results on almost periodic
trajectory and its local attractivity of shunting inhibitory
cellular neural networks (CNNs) with delays. In \cite{A2}, authors
proved existence and attractivty of almost periodic solution for
CNNs with distributed delays and variable coefficients.

In this paper,  we will investigate the existence of almost
periodic solution of (\ref{periodic}) and its global stability.
All the conclusions and methods used apply to discussion of
periodic solution when $d_{i}$, $a_{ij}$, $b_{ij}$ etc. are
periodic functions or equilibrium point when $d_{i}$, $a_{ij}$,
$b_{ij}$ etc. are constants.

\section{Some preliminaries}

\noindent{\bf Definition 1}\quad Class $H\{G_{1},\cdots,G_{n}\}$
of functions: Let ${G}=diag[{G}_{1},\cdots,{G}_{n}]$, where
${G}_{i}>0$, $i=1,\cdots,n$. $g(x)=(g_{1}(x),\cdots,g_{n}(x))^{T}$
is said to belong to $H\{G_{1},\cdots,G_{n}\}$, if the
non-decreasing functions $g_{i}(x)$, $i=1,\cdots,n$ satisfy
$\frac{|g_{i}(x+u)-g_{i}(x)|}{|u|}\le G_{i}$. Also,
$H_{1}\{F_{1},\cdots,F_{n}\}$ of functions: Let
${F}=diag[{F}_{1},\cdots,{F}_{n}]$, where ${F}_{i}>0$,
$i=1,\cdots,n$. $f(x)=(f_{1}(x),\cdots,f_{n}(x))^{T}$ is said to
belong to $H_{1}\{F_{1},\cdots,F_{n}\}$, if the functions
$f_{i}(x)$, $i=1,\cdots,n$ satisfy
$\frac{|f_{i}(x+u)-f_{i}(x)|}{|u|}\le F_{i}$.

\noindent{\bf Definition 2}\quad Let $x(t): R\rightarrow R^{n}$ be
continuous, $x(t)$ is said to be almost periodic (Bohr) on $R$ if
for any $\epsilon>0$, it is possible to find a real number
$l=l(\epsilon)>0$, for any interval with length $l(\epsilon)$,
there exists a number $\omega=\omega(\epsilon)$ in this interval
such that $||x(t+\omega)-x(t)||<\epsilon$, for $\forall t\in R$.

Throughout the paper, we use following norm
\begin{eqnarray*}
\|u(t)\|_{\{\xi,\infty\}}=\max_{i=1,\cdots,n}|\xi_{i}^{-1}u_{i}(t)|
\end{eqnarray*}

\section{Main Results}
In this section, we investigate the dynamical system
(\ref{almostperiodic}) with initial condition (\ref{ini}). In the
sequel, we make the following assumption.

\noindent{\bf Assumption A}
\begin{enumerate}
\item Activation functions $g(\cdot)\in H\{G_1,G_2,\cdots,G_n\}$,
$f(\cdot)\in H_{1}\{F_1,F_2,\cdots,F_n\}$ \item $d_{i}(t)$,
$a_{ij}(t)$, $\tau_{ij}(t)\geq 0$ and $I_{i}(t)$ are continuous,
and $d_{i}(t)\geq d_{i}>0$, for $i,j=1,2,\cdots,n$ \item For any
$s\in R$, $K_{ij}(t,s):t\rightarrow K_{ij}(t,s)$ is continuous,
and for any $t\in R$, $dK_{ij}(t,s): s\rightarrow K_{ij}(t,s)$ is
a Lebesgue-Stieljies measure, for all $i,j=1,2,\cdots,n$ \item For
any $\epsilon>0$, there exists $l=l(\epsilon)>0$, such that every
interval $[\alpha,\alpha+l]$  contains at least one number
$\omega$ for which
\begin{eqnarray}
&&|d_{i}(t+\omega)-d_{i}(t)|<\epsilon\nonumber\\
&&|a_{ij}(t+\omega)-a_{ij}(t)|<\epsilon\nonumber\\\
&&|I_{i}(t+\omega)-I_{i}(t)|<\epsilon\nonumber\\\
&&|\tau_{ij}(t+\omega)-\tau_{ij}(t)|<\epsilon\nonumber\\\
&&\int_{0}^{\infty}|dK_{ij}(t+\omega,s)-dK_{ij}(t,s)|<\epsilon\nonumber
\end{eqnarray}
hold for all $i,j=1.2.\cdots,n$ and $t\in R$.
\item $|dK_{ij}(t,s)|\le |dK_{ij}(s)|$, and for some
$\epsilon>0$ and $\int_{0}^{\infty}e^{\epsilon
s}|dK_{ij}(s)|<\infty$ .
\end{enumerate}
It can be seen that under item 4 in the assumption A,  $d_{i}(t)$,
$a_{ij}(t)$, $I_{i}(t)$ and $\tau_{ij}(t)$ are almost periodic
functions. Therefore, they are all bounded. We also denote
$|a_{ij}^{*}|=\sup_{\{t\in R \}}|a_{ij}(t)|<\infty,$
$|b_{ij}^{*}|=\sup_{\{t\in R
\}}\int_{0}^{\infty}|dK_{ij}(t,s)|<\infty,$
$|I_{i}^{*}|=\sup_{\{t\in R  \}}|I_{i}(t)|<\infty,$
$\tau^{\star}_{ij}=\sup_{\{t\in R \}}\tau_{ij}(t)<\infty$,
$i,j=1,\cdots,n$.

\begin{Lemma}\quad
Suppose that all items in Assumption A are satisfied. If there
exist $\xi_{i}>0$, $i=1,\cdots,n$, such that
\begin{equation}
-d_{i}(t)\xi_{i}+\sum\limits_{j=1
}^{n}|a_{ij}(t)|G_{j}\xi_{j}+\sum\limits_{j=1}^{n}F_{j}\xi_{j}\int_{0}^{\infty}|dK_{ij}(t,s)|
<-\eta<0
\end{equation}
hold for all $t>0$. In particular, if
\begin{equation}
-d_{i}\xi_{i}+\sum\limits_{j=1
}^{n}|a_{ij}^{*}|G_{j}\xi_{j}
+\sum\limits_{j=1}^{n}F_{j}\xi_{j}\int_{0}^{\infty}|dK_{ij}(s)|
<0
\end{equation}
Then any solution $u(t)$ of the system (\ref{almostperiodic})  is
bounded.
\end{Lemma}
{\bf Proof:}\quad Define $M(t)=\max\limits_{s\le
t}\|u(s)\|_{\{\xi,\infty\}}$. It is obvious that
$\|u(t)\|_{\{\xi,\infty\}}\le M(t)$, and $M(t)$ is non-decreasing.
We will prove that $M(t)\le \max\{M(0),\frac{2}{\eta}\hat{I}\}$,
where
\begin{eqnarray*}
\hat{I}&=&\max\limits_{i}\bigg\{|I^{*}_{i}|
+\sum\limits_{j=1}^{n}\bigg[|a_{ij}^{*}||g_{j}(0)|
+|b^{*}_{ij}||f_{j}(0)|\bigg]\bigg\}
\end{eqnarray*}

For any fixed $t_{0}$, there are two possible cases.

Case 1.
\begin{equation}
\| u(t_{0})\|_{\{\xi,\infty\}}<M(t_{0})=\max\limits_{s\le
t_{0}}\|u(s)\|_{\{\xi,\infty\}}
\end{equation}
In this case, in a small neighborhood $(t_{0},t_{0}+\delta)$ of
$t_{0}$, $\| u(t)\|_{\{\xi,\infty\}}<M(t_{0})$, and $M(t)=
M(t_{0})$

Case 2.
\begin{equation}
\| u(t_{0})\|_{\{\xi,\infty\}}=M(t_{0})=\max\limits_{s\le
t_{0}}\|u(s)\|_{\{\xi,\infty\}}
\end{equation}
In this case, let $i_{t_{0}}$ be such an index that
\begin{equation}
\xi_{i_{t_{0}}}^{-1}|u_{i_{t_{0}}}(t_{0})|=\|u(t_{0})\|_{\{\xi,\infty\}}
\end{equation}
Then noticing
\begin{eqnarray*}
|g_{j}(s)|\le G_{j}|s|+|g_{j}(0)|\quad |f_{j}(s)|\le
F_{j}|s|+|f_{j}(0)|\quad for~j=1,\cdots,n~s\in R
\end{eqnarray*}
we have
\begin{eqnarray}
&&\bigg\{\frac{d}{dt}|u_{i_{t_{0}}}(t)|\bigg\}_{t=t_{0}}\nonumber\\
&=&sign(u_{i_{t_{0}}}(t_{0}))\bigg[-d_{i_{t_{0}}}(t_{0})u_{i_{t_{0}}}(t_{0})
+\sum\limits_{j=1}^{n}a_{i_{t_{0}}j}(t_{0})g_{j}
(u_{j}(t_{0}))\nonumber\\
&+&\sum\limits_{j=1}^{n}
\int_{0}^{\infty}f_{j}(u_{j}(t_{0}-\tau_{i_{t_{0}}j}(t_{0})-s))
dK_{i_{t_{0}}j}(t_{0},s)
+I_{i_{t_{0}}}(t_{0})\bigg]\nonumber\\
&\le&-d_{i_{t_{0}}}(t_{0})|u_{i_{t_{0}}}(t_{0})|
\xi_{i_{t_{0}}}^{-1}\xi_{i_{t_{0}}}
+\sum\limits_{j=1
}^{n}|a_{i_{t_{0}}j}(t_{0})|G_{j}|u_{j}(t_{0})|
\xi_{j}^{-1}\xi_{j}\nonumber\\
&+&\sum\limits_{j=1}^{n}F_{j}\xi_{j}
\int_{0}^{\infty}|u_{j}(t_{0}-\tau_{i_{t_{0}}j}(t_{0})-s)|\xi_{j}^{-1}|
dK_{i_{t_{0}}j}(t_{0},s)|+|I_{i_{t_{0}}}(t_{0})|\nonumber\\
&&+\sum\limits_{j=1}^{n}|a_{i_{0}j}(t)||g_{j}(0)|
+\int_{0}^{\infty}|dk_{i_{0}j}(t,s)||f_{j}(0)|\nonumber\\
&\le&\bigg[-d_{i_{t_{0}}}(t_{0})\xi_{i_{t_{0}}} +\sum\limits_{j=1
}^{n}|a^{\star}_{i_{t_{0}}j}|G_{j}\xi_{j}\nonumber\\
&+&\sum\limits_{j=1}^{n}
F_{j}\xi_{j}\int_{0}^{\infty}|dK_{i_{t_{0}}j}(t_{0},s)|\bigg]
\|u(t_{0})\|_{\{\xi,\infty\}}+\hat{I}\nonumber\\
&\le&-\eta\|u(t_{0})\|_{\{\xi,\infty\}}+\hat{I}= -\eta M(t_{0})
+\hat{I}
\end{eqnarray}
Thus, if $M(t_{0})\geq \frac{2}{\eta}|\hat{I}|$, then $M(t)$ is
decreasing in a small neighborhood $(t_{0},t_{0}+\delta_{1})$ of
$t_{0}$. On the other hand, if $M(t_{0})<
\frac{2}{\eta}|\hat{I}|$, then in a small neighborhood
$(t_{0},t_{0}+\delta_{2})$ of $t_{0}$,
$\|u(t)\|_{\{\xi,\infty\}}<\frac{2}{\eta}|\hat{I}|$. Therefore,
$M(t)\le \max\{M(t_{0}),\frac{2}{\eta}|\hat{I}|\}$ in a small
neighborhood $(t_{0},t_{0}+\delta)$ of $t_{0}$, where
$\delta=min\{\delta_{1},\delta_{2})$.

In either case, we have $M(t)\le
\max\{M(t_{0}),\frac{2}{\eta}|\hat{I}|\}$ in a small neighborhood
$(t_{0},t_{0}+\delta)$ of $t_{0}$.

Starting from $t=0$, process this argument at every $t\geq 0$, we
conclude that if $M(0)>\frac{2}{\eta}|\hat{I}|$, then $M(t)=M(0)$
for all $t>0$. Instead, if $M(0)\le\frac{2}{\eta}|\hat{I}|$, then
$M(t)\le \frac{2}{\eta}|\hat{I}|$ for all $t>0$. Therefore,
$M(t)\le \max\{M(0),\frac{2}{\eta}|\hat{I}|\}$ for all $t>0$,
which proves that
 $u(t)$ is bounded. Lemma 1 is proved.

\begin{Lemma}\quad
Suppose that all items in Assumption A are satisfied. If there
exist $\xi_{i}>0$, $i=1,2,\cdots,n$, and $\beta>0$, such that fo
all $t>0$, there hold
\begin{eqnarray}
-d_{i}(t)\xi_{i}+\sum\limits_{j=1
i}^{n}|a_{ij}(t)|G_{j}\xi_{j}+\sum\limits_{j=1}^{n}F_{j}\xi_{j}e^{\beta
\tau^{\star}_{ij}} \int_{0}^{\infty} e^{\beta
s}|dK_{ij}(t,s)|<-\eta<0
\end{eqnarray}
Then for any $\epsilon>0$, there exist $T>0$ and
$l=l(\epsilon)>0$, such that every interval $[\alpha,\alpha+l]$
contains at least one number $\omega$ for which the solution
$u(t)$ of system (\ref{almostperiodic}) satisfies
\begin{equation}
\|u(t+\omega)-u(t)\|_{\{\xi,\infty\}}\le\epsilon\quad for ~all~t>T
\end{equation}
\end{Lemma}
{\bf Proof:}\quad 
Define
\begin{eqnarray*}
&&\epsilon_{i}(\omega,t)=-[d_{i}(t+\omega)-d_{i}(t)]
u_{i}(t+\omega)+\sum\limits_{j=1}^{n}[a_{ij}(t+\omega)-a_{ij}(t)]
g_{j}(u_{j}(t+\omega))\nonumber\\
&&+\sum\limits_{j=1}^{n}\int_{0}^{\infty}[f_{j}(u_{j}(t-\tau_{ij}(t+\omega)+\omega-s))
-f_{j}(u_{j}(t-\tau_{ij}(t)+\omega-s))]dK_{ij}(t+\omega,s)\nonumber\\
&&+\sum\limits_{j=1}^{n}\int_{0}^{\infty}f_{j}(u_{j}(t-\tau_{ij}(t)+\omega-s))d[K_{ij}(t+\omega,s)
-K_{ij}(t,s)]+[I_{i}(t+\omega)-I_{i}(t)]
\end{eqnarray*}
By Lemma 1, $u(t)$ is bounded. Thus, the right side of
(\ref{almostperiodic}) is also bounded, which implies that $u(t)$
is uniformly continuous. Therefore, by the assumption A, for any
$\epsilon>0$, there exists $l=l(\epsilon)>0$ such that every
interval $[\alpha,\alpha+l]$, $\alpha\in R$, contains an $\omega$
for which
\begin{equation}
 |\epsilon_{i}(\omega,t)|\le\frac{1}{2}\eta\epsilon\quad for~all~t\in
 R~,~i=1,2,\cdots,n
\end{equation}

Denote $x_{i}(t)=u_{i}(t+\omega)-u_{i}(t)$. We have
\begin{eqnarray}
\frac{dx_{i}(t)}{dt}&=&-d_{i}(t)x_{i}(t)+\sum\limits_{j=1}^{n}a_{ij}(t)
[g_{j}(u_{j}(t+\omega))
-g_{j}(u_{j}(t))]\nonumber\\
&&
+\sum\limits_{j=1}^{n}\int_{0}^{\infty}[f_{j}(u_{j}(t+\omega-\tau_{ij}(t)-s))-f_{j}(u_{j}
(t-\tau_{ij}(t)-s))]dK_{ij}(t,s)\nonumber\\
&&+\epsilon_{i}(\omega,t)
\end{eqnarray}
Let $i_{t}$ be such an index that
\begin{equation}
\xi_{i_t}^{-1}|x_{i_t}(t)|=\|x(t)\|_{\{\xi,\infty\}}
\end{equation}
Differentiate $e^{\beta s}|x_{i_t}(s)|$, we have
\begin{eqnarray}
&&\frac{d}{d s}\bigg\{e^{\beta
s}|x_{i_t}(s)|\bigg\}\bigg|_{s=t}=\beta e^{\beta
t}|x_{i_t}(t)|+e^{\beta
t}sign(x_{i_t}(t))\bigg\{-d_{i_t}(t)x_{i_t}(t)
\nonumber\\
&&+\sum\limits_{j=1}^{n}a_{i_{t}j}(t)\bigg[
g_{j}(u_{j}(t+\omega))-g_{j}(u_{j}(t))\bigg]\nonumber\\
&&+\sum\limits_{j=1}^{n}\int_{0}^{\infty}\bigg[
f_{j}(u_{j}(t+\omega-\tau_{i_{t}j}(t)-s))-f_{j}(u_{j}(t-\tau_{i_{t}j}(t)-s))\bigg]dK_{i_{t}j}(t,s)
+\epsilon_{i_{t}}(\omega,t)\bigg\}\nonumber\\
&&\le e^{\beta t}\bigg\{-[d_{i_t}(t)-\beta
]|x_{i_t}(t)|\xi_{i_t}^{-1} \xi_{i_t}+\sum\limits_{j=1
}^{n}|a_{i_{t}j}(t)|G_{j}|x_{j}(t)|\xi_{j}^{-1}\xi_{j}\nonumber\\
&&+\sum\limits_{j=1}^{n}F_{j}\xi_{j}
\int_{0}^{\infty}|x_{j}(t-\tau_{i_{t}j}(t)-s)|\xi_{j}^{-1}e^{-\beta(\tau_{i_{t}j}(t)+
s)}e^{\beta
(s+\tau^{\star}_{ij})}|dK_{i_{t}j}(t,s)|\bigg\}+\frac{1}{2}\eta\epsilon
e^{\beta t}
\end{eqnarray}

Similar to the proof of Lemma 1, let
\begin{eqnarray}
\Psi(t)=\max\limits_{s\le t }\bigg\{e^{\beta
s}\|x(s)\|_{\{\xi,\infty\}}\bigg\}
\end{eqnarray}

If there is such a point $t_0>0$ that $\Psi(t_0)=e^{\beta
t_0}\|x(t_0)\|_{\{\xi,\infty\}}$. Then we have
\begin{equation}
\frac{d}{dt}\bigg\{e^{\beta
t}|x_{i_{t_{0}}}(t)|\bigg\}_{t=t_0}<-\eta \Psi(t_0)+\eta\epsilon
e^{\beta t_0}
\end{equation}
In addition, if $\Psi(t_0)\ge \epsilon e^{\beta t_0}$, then
$\Psi(t)$ is decreasing in a small neighborhood
$(t_{0},t_{0}+\delta)$ of $t_0$. On the other hand, if $\Psi(t_0)<
\epsilon e^{\beta t_0}$, then in a small neighborhood
$(t_{0},t_{0}+\delta)$ of $t_{0}$, $e^{\beta
t}\|x(t)\|_{\{\xi,\infty\}}<\epsilon e^{\beta t_0}$, and
$\Psi(t)<\max\{\Psi(t_0),\epsilon e^{\beta t_0}\}$. By the same
reasonings used in the proof of Lemma 1, for all $t>t_0$, we have
$e^{\beta t }\|x(t)\|_{\{\xi,\infty\}}\le\max\{\Psi(t_0),\epsilon
e^{\beta t }\}$. Therefore, there exists $t_1>0$, for all $t>t_1$,
$\|x(t)\|_{\{\xi,\infty\}}\le\epsilon$.

Instead, if for all $t>0$, we have $\Psi(t)>e^{\beta
t}\|x(t)\|_{\{\xi,\infty\}}$, then $\Psi(t)=\Psi(0)$ is a
constant, and $e^{\beta t}\|x(t)\|_{\{\xi,\infty\}}\le \Psi (t)=
\Psi (0)$. Hence, there exists $t_{2}>0$, for all $t>t_{2}$, such
that $\|x(t)\|_{\{\xi,\infty\}}\le\epsilon$.

In summary, there must exist $T>0$ such that
$\|x(t)\|_{\{\xi,\infty\}}\le\epsilon$ holds for all $t>T$. Lemma
2 is proved.

\begin{Theorem}\quad Suppose that all items in Assumption A are satisfied. If
there exist $\xi_{i}>0$, $i=1,2,\cdots,n$, and $\beta>0$ such that
\begin{eqnarray}
&-&[d_{i}(t)-\beta]\xi_{i}+\sum\limits_{j=1
}^{n}|a_{ij}(t)|G_{j}\xi_{j}
+\sum\limits_{j=1}^{n}F_{j}\xi_{j}e^{\beta \tau^{\star}_{ij}}
\int_{0}^{\infty} e^{\beta s}|dK_{ij}(t,s)|< 0 \label{Th1}
\end{eqnarray}
hold for all $t>0$. In particular,
\begin{eqnarray}
&-&[d_{i}-\beta]\xi_{i}+\sum\limits_{j=1
}^{n}|a^{\star}_{ij}|G_{j}\xi_{j}
+\sum\limits_{j=1}^{n}F_{j}\xi_{j}e^{\beta \tau^{\star}_{ij}}
\int_{0}^{\infty} e^{\beta s}|dK_{ij}(s)|< 0 \label{Th11}
\end{eqnarray}
Then the dynamical system (\ref{almostperiodic})
has a unique almost periodic solution $v(t)=[v_{1}(t),
v_{2}(t),\newline \cdots,v_{n}(t)]^{T}$, and for any solution
$u(t)=[u_{1}(t),u_{2}(t),\cdots,u_{n}(t)]^{T}$ of
(\ref{almostperiodic}), there holds
\begin{equation}
\|u(t)-v(t)\|=O(e^{-\beta t})
\end{equation}
\end{Theorem}
{\bf Proof:}\quad As did before, define
\begin{eqnarray}
&&\epsilon_{i,k}(t)=-[d_{i}(t+t_{k})-d_{i}(t)]
u_{i}(t+t_{k})+\sum\limits_{j=1}^{n}[a_{ij}(t+t_{k})-a_{ij}(t)]
g_{j}(u_{j}(t+t_{k}))\nonumber\\
&&+\sum\limits_{j=1}^{n}\int_{0}^{\infty}[f_{j}(u_{j}(t-\tau_{ij}(t+t_{k})+t_{k}-s))
-f_{j}(u_{j}(t-\tau_{ij}(t)+t_{k}-s))]dK_{ij}(t+t_{k},s)\nonumber\\
&&+\sum\limits_{j=1}^{n}\int_{0}^{\infty}f_{j}(u_{j}(t+t_{k}-\tau_{ij}(t)-s))
d[K_{ij}(t+t_{k},s)
-K_{ij}(t,s)]+[I_{i}(t+t_{k})-I_{i}(t)]\nonumber\\
\end{eqnarray}
where $t_{k}$ is any sequence of real numbers. From the Assumption
A and the boundedness of $u(t)$, we can select a sequence
$\{t_{k}\}\rightarrow\infty$ such that
\begin{equation}
|\epsilon_{i,k}(t)|\le\frac{1}{k}\quad for ~all~i,t
\end{equation}
Because $\{u(t+t_{k})\}_{k=1}^{\infty}$ are uniformly bounded and
equiuniformly continuous. By Arzala-Ascoli Lemma and diagonal
selection principle, we can select a subsequence $t_{k_{j}}$ of
$t_{k}$, such that $u(t+t_{k_{j}})$ (for convenience, we still
denote by $u(t+t_{k})$) uniformly converges to a continuous
function $v(t)=[v_{1}(t),v_{2}(t),\cdots,v_{n}(t)]^T$ on any
compact set of $R$.

Now, we  prove $v(t)$ is a solution of system
(\ref{almostperiodic}). In fact, by Lebesgue dominant convergence
theorem, for any $t>0$ and $\delta t\in R$, we have
\begin{eqnarray}
&&v_{i}(t+\delta
t)-v_{i}(t)=\lim\limits_{k\rightarrow\infty}\bigg[u_{i}(t+\delta
t+t_{k})-u_{i}(t+t_{k}) )\bigg]\nonumber\\
&=&\lim\limits_{k\rightarrow\infty}\int_{t}^{t+\delta
t}\bigg\{-d_{i}(\sigma+t_{k})u_{i}(\sigma+t_{k})+\sum\limits_{j=1}^{n}
a_{ij}(\sigma+t_{k})
g_{j}(u_{j}(\sigma+t_{k}))\nonumber\\
&&+\sum\limits_{j=1}^{n}\int_{0}^{\infty}f_{j}(u_{j}(\sigma+t_{k}-\tau_{ij}(\sigma+t_{k})-s))
dK_{ij}(\sigma+t_{k},s)+I_{i}(\sigma+t_{k})\bigg\}d\sigma\nonumber\\
&=&\int_{t}^{t+\delta
t}\bigg\{-d_{i}(\sigma)v_{i}(\sigma)+\sum\limits_{j=1}^{n}a_{ij}(\sigma)
g_{j}(v_{j}(\sigma))\nonumber\\
&&+\sum\limits_{j=1}^{n}\int_{0}^{\infty}f_{j}(v_{j}(\sigma-\tau_{ij}(\sigma)-s))
dK_{ij}(\sigma,s)
+I_{i}(\sigma)\bigg\}d\sigma\nonumber\\
&&+\lim\limits_{k\rightarrow\infty}\int_{t}^{t+\delta t
}\epsilon_{i,k}(s)d\sigma\nonumber\\
&=&\int_{t}^{t+\delta
t}\bigg\{-d_{i}(\sigma)v_{i}(\sigma)+\sum\limits_{j=1}^{n}a_{ij}(\sigma)
g_{j}(v_{j}(\sigma))\nonumber\\
&&+\sum\limits_{j=1}^{n}\int_{0}^{\infty}f_{j}(v_{j}(\sigma-\tau_{ij}(\sigma)-s))
dK_{ij}(\sigma,s) +I_{i}(\sigma)\bigg\}d\sigma
\end{eqnarray}
which implies
\begin{equation}
\frac{dv_{i}}{dt}=-d_{i}(t)v_{i}(t)+\sum\limits_{j=1}^{n}a_{ij}(t)
g_{j}(u_{j}(t))
+\int_{0}^{\infty}f_{j}(u_{j}(t-\tau_{ij}(t)-s))dK_{ij}(t,s)+I_{i}(t)
\end{equation}
Therefore, $v(t)$ is a solution of system (\ref{almostperiodic}).

Secondly, we prove that $v(t)$ is an almost periodic function. By
Lemma 2, for any $\epsilon>0$, there exist $T>0$ and
$l=l(\epsilon)>0$, such that every interval $[\alpha,\alpha+l]$
contains at least one number $\omega$ for which
\begin{equation}
|u_{i}(t+\omega)-u_{i}(t)|\le\epsilon\quad for~ all ~t>T
\end{equation}
Then we can find a sufficient large $K\in N$ such that for any
$k>K$
\begin{equation}
|u_{i}(t+t_{k}+\omega)-u_{i}(t+t_{k})|\le\epsilon\quad for~all~
t>0
\end{equation}
holds. Let $k\rightarrow\infty$, we have
\begin{equation}
|v_{i}(t+\omega)-v_{i}(t)|\le\epsilon\quad for~all~t>0
\end{equation}
In other words, $v(t)$ is an almost periodic function.

Finally, we prove every solution $u(t)$ of (\ref{almostperiodic})
converges to $v(t)$ exponentially with rate $\beta$.

Denote $y(t)=u(t)-v(t)$, we have
\begin{eqnarray}
\frac{dy_{i}(t)}{dt}&=&-d_{i}(t)y_{i}(t)+\sum\limits_{j=1}^{n}a_{ij}(t)
[g_{j}(u_{j}(t))
-g_{j}(v_{j}(t))]\nonumber\\
&&+\sum\limits_{j=1}^{n}\int_{0}^{\infty}[f_{j}(u_{j}(t-\tau_{ij}(t)-s))
-f_{j}(v_{j}(t-\tau_{ij}(t)-s))] dK_{ij}(t,s)
\end{eqnarray}

Let $i_{t}$ be such an index that
\begin{equation}
|y_{i_t}(t)|=\xi_{i_t}\|y(t)\|_{\{\xi,\infty\}}
\end{equation}
and differentiate $e^{\beta s}|y_{i_t}(s)|$, we have
\begin{eqnarray}
&&\frac{d}{ds}\bigg\{e^{\beta
s}|y_{i_t}(s)|\bigg\}\bigg|_{s=t}=\beta e^{\beta
t}|y_{i_t}(t)|+e^{\beta
t}sign(y_{i_t}(t))\bigg\{-d_{i_t}(t)y_{i_t}(t)\nonumber\\
&& +\sum\limits_{j=1
}^{n}a_{i_{t}j}(t)[g_{j}(u_{j}(t))-g_{j}(v_{j}(t))]\nonumber\\
&&+\sum\limits_{j=1}^{n}
[f_{j}(u_{j}(t-\tau_{i_{t}j}(t)-s))-f_{j}(v_{j}(t-\tau_{i_{t}j}(t)-s))]dK_{i_{t}j}(t,s)
\bigg\}\nonumber\\
&\le&e^{\beta t}\bigg\{-[d_{i_t}-\beta]|y_{i_t}(t)| \xi_{i_t}^{-1}
\xi_{i_t}+\sum\limits_{j=1
}^{n}|a_{i_{t}j}(t)|G_{j}|y_{j}(t)|\xi_{j}^{-1}\xi_{j}\nonumber\\
&&+\sum\limits_{j=1}^{n}
F_{j}\xi_{j}\int_{0}^{\infty}|y_{j}(t-\tau_{i_{t}j}(t)-s)|\xi_{j}^{-1}e^{-\beta
(s+\tau_{i_{t}j}(t)) }e^{\beta
(s+\tau^{\star}_{i_{t}j})}|dK_{i_{t}j}(t,s)|\bigg\}
\label{convergence}
\end{eqnarray}

Define $\Delta(t)=\max\limits_{s\le t }\bigg\{e^{\beta
s}\|y(s)\|_{\{\xi,\infty\}}\bigg\}$. If for some $t_{0}$,
\begin{eqnarray}
\Delta(t_0)=e^{\beta t_0}\|y(t_0)\|_{\{\xi,\infty\}}
\end{eqnarray}
Then, by (\ref{convergence}), we obtain
\begin{equation}
\frac{d}{dt}\{e^{\beta
t}|y_{i_0}(t)|\}_{t=t_0}\le-\eta\Delta(t_0)\le 0
\end{equation}
and $\Delta(t)=\Delta(t_0)$ in a small neighborhood
$(t_{0},t_{0}+\delta)$ of $t_{0}$. On the other hand, if
\begin{eqnarray}
e^{\beta t_0}\|y(t_0)\|_{\{\xi,\infty\}}<\Delta(t_0)
\end{eqnarray}
Then in a small neighborhood of $t_{0}$,
\begin{eqnarray}
e^{\beta t}\|y(t)\|_{\{\xi,\infty\}}<\Delta(t_0)
\end{eqnarray}
Therefore, $\Delta(t)=\Delta(0)$ for all $t\ge 0$, and
\begin{equation}
\|y(t)\|_{\{\xi,\infty\}}\le \Delta(0)e^{-\beta t}
\end{equation}
Namely,
\begin{equation}
\|u(t)-v(t)\|_{\{\xi,\infty\}}\le \Delta(0)e^{-\beta t}
\end{equation}
Theorem 1 is proved.

{\bf Remark 1}\quad From the properties of continuous and almost
periodic functions, i.e., Bochner's Theorem, see (pp.4-7,
\cite{Levitan}) it can be seen that if $d_{i}(t)$, $a_{ij}(t)$,
$b_{ij}(t)$, $\tau_{ij}(t)$ and $I_{i}(t)$ are all continuous and
almost periodic, then the item $4$ in Assumption (A) are
satisfied.

If $dK_{ij}(t,0)=b_{ij}(t)$, and $dK_{ij}(t,s)=0$, for $s\ne 0$,
then (\ref{almostperiodic}) reduces to
\begin{eqnarray}
\frac{du_i}{dt}=-d_{i}(t)u_{i}(t)+\sum_{j=1}^{n}a_{ij}(t)g_j(u_j(t))
+\sum_{j=1}^{n}b_{ij}(t)f_{j}(u_{j}(t-\tau_{ij}(t)))+I_i(t)
\label{Co1}
\end{eqnarray}
Therefore, we have following

{\bf Corollary 1}\quad
 Suppose $g(x)=(g_{1}(x),\cdots,g_{n}(x))^{T}\in$
$H\{G_{1},\cdots,G_{n}\}$ and $f(x)=
(f_{1}(x),\cdots,f_{n}(x))^{T}\in H\{F_{1},\cdots,F_{n}\}$.
$d_{i}(t)$, $a_{ij}(t)$, $b_{ij}(t)$, $\tau_{ij}(t)\geq 0$ and
$I_{i}(t)$, $i,j=1,\cdots,n$,  are continuous almost periodic
functions. If there are positive constants
$\xi_{1},\cdots,\xi_{n}$ and $\beta$ such that
\begin{eqnarray}
\xi_{i}[-d_{i}(t)+\beta] +\sum_{j=1 }^{n}\xi_{j}|a_{ij}(t)|G_{j}
+\sum_{j=1}^{n}|b_{ij}(t)| F_{j}e^{\beta \tau^{\star}_{ij}}< 0,
\quad i=1,\cdots,n,
\end{eqnarray}
hold for all $t>0$. In particular, if
\begin{eqnarray}
\xi_{i}[-d_{i}+\beta] +\sum_{j=1 }^{n}\xi_{j}|a_{ij}^{*}|G_{j}
+\sum_{j=1}^{n}|b_{ij}^{*}| F_{j}e^{\beta \tau^{\star}_{ij}}< 0,
\quad i=1,\cdots,n,
\end{eqnarray}Then the dynamical system (\ref{Co1})
 has a unique almost periodic
solution  $v(t)$ and, for any solution $u(t)$ of (\ref{Co1}),
there holds
\begin{equation}
||u(t)-v(t)||=O(e^{-\beta t})
\end{equation}

Instead, if $dK_{ij}(t,s)=b_{ij}(t)k_{ij}(s)ds$, then we have
 the following distributed delayed dynamical systems:
\begin{eqnarray}
&&\frac{du_{i}(t)}{dt}=-d_{i}(t)u_{i}(t)+\sum_{j=1}^{n}a_{ij}(t)g_{j}(u_{j}(t))\nonumber\\
&+&\sum_{j=1}^{n}b_{ij}(t)\int_{0}^{\infty}k_{ij}(s)f_{j}(u_{j}(t-\tau_{ij}(t)-s))ds
+I_{i}(t), \quad i=1,\cdots,n \label{Co2}
\end{eqnarray}

\noindent{\bf Corollary 2}\quad Suppose
$g(x)=(g_{1}(x),\cdots,g_{n}(x))^{T}\in H\{G_{1},\cdots,G_{n}\}$
and $f(x)= (f_{1}(x),\cdots,f_{n}(x))^{T}\in
H\{F_{1},\cdots,F_{n}\}$,
 $d_{i}(t)$, $a_{ij}(t)$,
$b_{ij}(t)$, $k_{ij}(t)$, $\tau_{ij}(t)\geq 0$ and $I_{i}(t)$,
$i,j=1,\cdots,n$,
 are continuous almost periodic functions and there
exists $\beta>0$ such that
\begin{eqnarray}
\int_{0}^{\infty} e^{\beta s}|k_{ij}(s)|ds<\infty
\end{eqnarray}
If
there are positive constants $\xi_{1},\cdots,\xi_{n}$ such that
\begin{eqnarray}
\xi_{i}[-d_{i}+\beta] +\sum_{j=1 }^{n}\xi_{j}|a_{ij}(t)|G_{j}
+\sum_{j=1}^{n}|b_{ij}(t)| F_{j}e^{\beta \tau^{\star}_{ij}}
\int_{0}^{\infty} e^{\beta s}|k_{ij}(s)|ds\le 0,
\end{eqnarray}
hold for all $t>0$. In particular, if
\begin{eqnarray}
\xi_{i}[-d_{i}+\beta] +\sum_{j=1 }^{n}\xi_{j}|a_{ij}^{*}|G_{j}
+\sum_{j=1}^{n}|b_{ij}^{*}| F_{j}e^{\beta \tau^{\star}_{ij}}
\int_{0}^{\infty} e^{\beta s}|k_{ij}(s)|ds\le 0,
\end{eqnarray}
Then the dynamical system (\ref{Co2})
 has a unique almost periodic
solution  $v(t)$ and, for any solution $u(t)$ of (\ref{Co2}), we
have
\begin{equation}
||u(t)-v(t)||=O(e^{-\beta t})
\end{equation}

{\bf Remark 2}\quad In \cite{Huang,A1,A2}, authors discussed
almost periodic solution for CNNs with delays (constant delays
$\tau_{ij}$ in \cite{Huang}, $\tau_{ij}(t)=\tau(t)$ in \cite{A1},
and distributed delay with $\tau_{ij}(t)=0$ in \cite{A2}) and its
local stability. Our model unifies models with time-varying delays
$\tau_{ij}(t)$ and distributed delays and even both, and includes
all the systems in these three papers as special cases. The
stability is global. Moreover, in these three papers, all the
reasonings heavily depend on the Theorem 1.93 in \cite{He}. Thus,
many unreasonable restrictions are imposed on the activation
functions and the interconnection weights (for example,
interconnection weights must satisfy $|a_{ij}(t)|>c$ and
$|b_{ij}(t)|>c$ for some positive constant and many others). In
this paper, we propose a totally different approach to prove
existence of almost periodic solution and its global stability
without citing any complicated theory. As we show, all the
redundant restrictions needed in these papers are removed.

\end{document}